\documentclass[11pt]{amsart}
\usepackage[active]{srcltx}
\usepackage{calc,amssymb,amsthm,amsmath,amscd, eucal,ulem}
\usepackage{alltt}
\usepackage[left=1.35in,top=1.25in,right=1.35in,bottom=1.25in]{geometry}
\RequirePackage[dvipsnames,usenames]{color}

\usepackage{tikz}

\usepackage{fullpage}
\usepackage{setspace}
\usepackage{hyperref}
\usepackage{enumerate}

\usepackage[all,cmtip]{xy}

\newcommand{\m}{\mathfrak{m} }

\providecommand{\e}{{\mathcal E}}
\providecommand{\D}{{\mathcal D}}
\newcommand{\Z}{\mathbb{Z} }
\newcommand{\N}{\mathbb{N} }

\newcommand{\TT}{\mathcal{T}}

\newcommand{\rt}{\rightarrow}

\newcommand{\chars}{\operatorname{char}}

\providecommand\Mod{\text{\rm Mod}}

\providecommand\Spec{\text{\rm Spec}}

\newcommand{\Hom}{\operatorname{Hom}}

\newcommand{\Ext}{\operatorname{Ext}}
\newcommand{\Tor}{\operatorname{Tor}}
\newcommand{\Tot}{\operatorname{Tot}}
\newcommand{\id}{\operatorname{id}}

\theoremstyle{plain}

\newtheorem{theorem}{Theorem}[section]
\newtheorem{corollary}[theorem]{Corollary}
\newtheorem{lemma}[theorem]{Lemma}
\newtheorem{proposition}[theorem]{Proposition}

\theoremstyle{definition}
\newtheorem{definition}[theorem]{Definition}
\newtheorem{s}[theorem]{}
\newtheorem{remark}[theorem]{Remark}

\newtheorem{remarks}[theorem]{Remarks}
\newtheorem*{example*}{Example}
\newtheorem*{claim*}{\it Claim}
\newtheorem*{note*}{\it Note}

\begin{document}
\title[On derived functors]{On derived functors of graded local cohomology modules-II}
\address{Department of Mathematics, Indian Institute of Technology Bombay, Powai, Mumbai 400 076, India}
\author{Tony J. Puthenpurakal}
\email{\href{mailto:tputhen@math.iitb.ac.in}{tputhen@math.iitb.ac.in}}
\author{Sudeshna Roy}
\email{\href{mailto:sudeshnaroy.11@gmail.com}{sudeshnaroy.11@gmail.com}}
\address{Department of Mathematics, Visvesvaraya National Institute of Technology, South Ambazari Road, Nagpur 440 010, India}
\author{Jyoti Singh}
\email{\href{mailto:jyotijagrati@gmail.com}{jyotijagrati@gmail.com}}
\date{\today}

\begin{abstract}

Let $R=K[X_1,\ldots, X_n]$ where $K$ is a field of characteristic zero, and let $A_n(K)$ be the $n^{th}$ Weyl algebra over $K$. We give  standard grading on $R$ and $A_n(K)$. Let $I$, $J$ be homogeneous ideals of $R$. Let $M = H^i_I(R)$ and $N = H^j_J(R)$ for some $i, j$. We show that $\Ext_{A_n(K)}^{\nu}(M,N)$ is concentrated in degree zero for all $\nu \geq 0$, i.e., $\Ext_{A_n(K)}^{\nu}(M,N)_l=0$ for $l \neq0$. This proves a conjecture stated in part I of this paper.
\end{abstract}

\maketitle
\section{Introduction}

\s {\it Setup}:\label{setup} Let $K$ be a field of characteristic zero and $R=K[X_1, \ldots,X_n]$ with  standard grading. Let $A_n(K)=K \langle X_1,\ldots,X_n, \partial_1, \ldots, \partial_n \rangle$  be the $n^{th}$ Weyl algebra over $K$. Consider $A_n(K)$ as graded with $ \deg X_i = 1$  and $\deg \partial_i = -1$ for $i = 1,\ldots, n$. Let $I, J$ be homogeneous ideals in $R$. Fix $i, j \geq 0$ and set $M = H^i_I(R)$ and $N = H^j_J(R)$.  In \cite[2.9]{Lyu-1}, Lyubeznik showed that the local cohomology module $H^i_I(R)$ is a 
graded left holonomic $A_n(K)$-module for $i \geq 0$.  

In \cite{TP6}, the first author proved that the de Rham cohomology $H^{\nu}(\partial, H_I^i(R))$ is concentrated in degree $-n$ for each $\nu \geq 0$. We have a graded isomorphism 
\begin{equation}\label{eq1}
H^{\nu}(\partial, H^i_I(R)) \cong \Tor_{n -\nu}^{A_n(K)}(R^r, H_I^i(R)),
\end{equation}
where $R^r$ is $R$ considered as a right $A_n(K)$-module by the isomorphism
$R^r = A_n(K)/(\partial)A_n(K)$. By $ ^lR$ we denote $R$ considered as a left 
$A_n(K)$-module  via the isomorphism $ ^lR \cong A_n(K)/A_n(K)(\partial)$. Let $L^{\sharp}$ denote the standard right $A_n(K)$-module associated to the left $A_n(K)$-module $L$ (see \ref{trans}). Then it is easy to check that $( ^lR)^\sharp \cong R^r$. Note that $H^0_{(0)}(R) =\ ^lR$. In \cite{TP5}, the first author along with the third author considerably generalized \eqref{eq1} and proved the following:

\begin{theorem}[with hypothesis as in \ref{setup}]\label{tor-con}\cite[Theorem 1.1]{TP5}
Fix $\nu \geq 0$. Then the graded $K$ vector space  $\Tor^{A_n(K)}_\nu(M^\sharp, N)$ 
is concentrated in degree $-n$, i.e., $\Tor^{A_n(K)}_\nu(M^\sharp, N)_j = 0$ for all $j \neq -n$.
\end{theorem}

In view of Theorem \ref{tor-con}, it is natural to investigate $\Ext^{\nu}_{A_n(K)}(M,N)$ which is a finite dimensional $K$-vector space for any $\nu \geq 0$, see \cite[2.7.15, 1.6.6]{BJ}. 

In \cite[Section 8]{TP5}, the first and third author conjectured that $\Ext_{A_n(K)}^{\nu}(M,N)$ is concentrated in degree zero for all $\nu \geq 0$. They gave several examples in support of this conjecture. The main result of this paper is an affirmative answer to this conjecture, i.e., we have

\begin{theorem}[with hypotheses as in \ref{setup}.] \label{conj}\cite[Conjecture 1.5]{TP5} Fix $\nu \geq 0$. The graded $K$-vector space $\Ext^\nu_{A_n(K)}(M, N)$ is concentrated in degree zero, i.e., $\Ext^\nu_{A_n(K)}(M, N)_l = 0$  for $l \neq 0$.
\end{theorem}

\s {\it Techniques used to prove Theorem \ref{conj}:}

Consider the Eulerian operator $\e_n =\sum_{i=1}^{n}X_i \partial_i$ and set $|w|=\deg w$ for any homogeneous element $w$ of a graded $A_n(K)$-module $W$. We say a graded $A_n(K)$-module $W$ is {\it Eulerain} if $\e_n w = |w|w$ for each homogeneous element $w$ of $W$. We say $W$ is {\it generalized Eulerian} if for each homogeneous element $w$ of $W$ there exists a positive integer $a$ depending on $w$ such that $(\e_n -|w|)^aw = 0$. 

In this paper we define the following. 

\begin{definition}\label{SGE}
We say a graded $A_n(K)$-module $W$ is {\it strongly generalized Eulerian} if for all homogeneous element $w$ there exists $a>0$ independent of $w$ such that 
\begin{equation}\label{strong-ge-def}
(\e_n - |w|)^aw = 0.
\end{equation}
\end{definition}

The notion of Eulerian modules was introduced by L. Ma and W. Zhang in \cite{MZ} over $\D$, the ring of K-linear differential operators on $R = K[X_1,\cdots,X_n]$ for both the cases when $K$ has characteristic zero and 
$p > 0$. But the class of Eulerian $\D$-modules is not closed under extensions (see \cite[3.5(1)]{MZ}). To rectify this, in \cite{TP6} the first author introduced the notion of generalized Eulerian $\D$-modules for the case of characteristic zero. {If $\chars K=0$, then $\D \cong A_n(K)$.} In this paper, we introduce  the notion of strongly generalized Eulerian $A_n(K)$-modules.

Recall that if a generalized Eulerian $A_n(K)$-module $M$ is finitely generated as an $R$-module, then $M = 0$ or $M = R^m$ for some $m \geq 1$, see \cite[Corollary 5.6]{TP5}. 
In Lemma \ref{ge-equiv-sge}, we show that for any finitely generated $A_n(K)$-module (e.g., holonomic $A_n(K)$-modules) these two notions generalized Eulerian and strongly generalized Eulerian are equivalent. In particular, we get the following results: 
\begin{enumerate}[\rm (a)]
	\item If $\TT$ is a graded {\it Lyubeznik functor} on ${}^* \Mod(R)$, then $\TT(R)$ is a strongly generalized Eulerian $A_n(K)$-module (as $\TT(R)$ is holonomic by \cite[2.2(d)]{Lyu-1}).
	\item If $M$ and $N$ are graded holonomic generalized Eulerian $A_n(K)$-modules, 
then $\Tor^R_{i}(M,N)$ is a holonomic strongly generalized Eulerian $A_n(K)$-module for all $i \geq 0$. \end{enumerate}

In fact we do not have an example of a generalized Eulerian module which is not strongly generalized Eulerian. 

In \cite[Proposition 3.9]{SZ}, N. Switala and W. Zhang showed that if $M$ is an Eulerian graded $A_n(K)$-module, then so is $D(M):={}^*\Hom_R(M, E)$, where $E=H^n_{(X_1, \ldots, X_n)}(R)$. From \cite[Proposition 2.13]{SZ} we have $D(M)\cong {}^*\Hom_K(M(-n), K)=M(-n)^\vee$. Inspired by their result, we prove the following:

\begin{theorem}\label{main}
Let $M^{\vee} :={}^* \Hom _K(M, K)$ denote the graded Matlis dual of $M$. If $M$ is a strongly generalized Eulerian $A_n(K)$-module, then so is $M(-n)^{\vee}$.
\end{theorem}

As an application of Theorem \ref{main} we get the following result:

\begin{corollary}\label{ext}
Let $N$ be strongly generalized Eulerian and $M$ be generalized Eulerian $A_n(K)$-modules. Then ${}^*\Ext^{\nu}_R(M,N)$ is generalized Eulerian for any $\nu \geq 0$.
\end{corollary}

We use Corollary \ref{ext} and a spectral sequence argument to prove our main result Theorem \ref{conj}.

\vspace{0.25cm}
The paper is structured as follows: In Section 2 we recall some definitions and basic facts used in the following sections. In Section 3 we discuss the properties of strongly generalized Eulerian modules and prove that for any finitely generated $A_n(K)$-module, the notions generalized Eulerian and strongly generalized Eulerian are equivalent. In Section 4 we prove Theorem \ref{main} and Corollary \ref{ext}. Section 5 is devoted to proving the main result Theorem \ref{conj}.


\section{Preliminaries}
In this section, we discuss a few preliminary result that we need.

\s {\bf Graded Lyubeznik functors}
Let $K$ be a field (not necessarily of characteristic zero)
and let $R=K[X_1, \ldots, X_n]$ be standard graded. We say $Y$ is homogeneous closed subset of $\Spec(R)$ if $Y=V(f_1, \ldots, f_s)$, where $f_i$'s are homogeneous polynomials in $R$. We say $Y$ is homogeneous locally closed subset of $\Spec(R)$ if $Y=Y''-Y'$ where $Y'\subset Y''$ are homogeneous closed subsets of $\Spec(R)$. Let $^*\Mod(R)$ be the category of graded $R$-modules. Then we have an exact sequence of functors on $^*\Mod(R)$, 
\begin{equation}\label{glyu}
\cdots \to H_{Y'}^i(-) \to H_{Y''}^i(-) \to H_Y^i(-) \to H_{Y'}^{i+1}(-) \to \cdots.
\end{equation}
\begin{definition}
	A {\it graded Lyubeznik functor} $\TT$ is a composite functor of the form $\TT=\TT_1 \circ \TT_2 \circ \cdots \circ\TT_m$ where each $\TT_j$ is either $H_{Y_j}^i(-)$ for some homogeneous locally closed subset $Y_j$ of $\Spec(R)$ or the kernel, image or cokernel of any arrow appearing in \eqref{glyu} with $Y=Y_j, ~Y'=Y'_j$ and $Y''=Y''_j$ such that $Y_j=Y''_j-Y'_j$ and $Y'_j \subset Y''_j$ are homogeneous closed subsets of $\Spec(R)$.
\end{definition}

\s Let $K$ be a field of characteristic zero and  let $R=K[X_1, \ldots, X_n]$  with standard grading. Let $A_n(K)$ be the $n^{th}$ Weyl algebra over $K$. Set $\deg \partial_i= -1$. So $A_n(K)$ is a graded ring.
The Euler operator, denoted by $\mathcal{E}_n$, is defined as $\e_n=\sum_{i=1}^{n}X_i \partial_i$.
Note that $\deg \mathcal{E}_n=0$. Let $M$ be a graded $A_n(K)$-module. If $z \in M$ is a homogeneous element, set $|z|= \deg z$.

\begin{definition}
 Let $M$ be a graded $A_n(K)$-module. We say $M$ is  \textit{Eulerian} if for any  homogeneous element $z\in M$, 
$$ \mathcal{E}_nz= |z|\cdot z.$$
\end{definition}

\begin{definition}
A graded $A_n(K)$-module $M$ is said to be \textit{generalized Eulerian} if for a homogeneous element $z$ of $M$ there exists a positive integer $a$ depending on $z$ such that
$$ (\mathcal{E}_n- {|z|})^a \cdot z = 0.$$
\end{definition}

\begin{lemma}\label{loca}
Let $M$ be a graded $A_n(K)$-module. Consider a homogeneous polynomial $f \in R$ and let $z$ be a homogeneous element of $M$. Then
$$ (\mathcal{E}_n- |z| + |f|)^n\frac{z}{f} = \frac{1}{f}(\mathcal{E}_n- |z|)^nz, ~~~~~ \text{for all}~~~~~ n\geqslant 1. $$
\end{lemma}

\begin{proof}
See Lemma 3.6 in \cite{TP5}.
\end{proof}

\s \label{ext-pre} 
Let $S= \bigoplus_{i \in \Z} S_i$ be a graded ring (not necessarily standard graded). A graded $S$-module $M$ is a $S$-module with a decomposition $M=\bigoplus_{i \in \Z} M_i$ such that $S_iM_j \subseteq M_{i+j}$ for all $i, j \in \Z$. We denote the category of graded $S$-modules by ${}^*\Mod(S)$. Here the objects are the graded $S$-modules and morphism $\phi: M \to N$ is a $S$-module homomorphism satisfying $\phi(M_i) \subseteq N_i$ for all $i \in \Z$.

Let $M, N$ be graded $S$-modules. A $S$-module homomorphism $\phi: M \to N$ is called homogeneous of degree $i$ if $M_n \subseteq N_{n+i}$ for all $i$. By $\Hom_i(M,N)$ we denote the group of homogeneous homomorphisms of degree $i$. Set $ ^*\Hom_S(M, N)=\bigoplus_{i \in \Z} \Hom_i(M,N)$. In general, $^*\Hom_S(M, N)$ is a $S$-submodule of $\Hom_S(M,N)$. For any $N \in {}^*\Mod(S)$ we define $ ^*\Ext^i_S(M,N)$ to be the $i^{th}$ right derived functor of $^*\Hom_S(M, N)$ in $ ^*\Mod(S)$. Thus if $\mathbb{I}^{\bullet}$ is a graded injective resolution of $N$, then $$ ^*\Ext^i_S(M,N) \cong H^i( ^*\Hom_S(M,\mathbb{I}^{\bullet})).$$ From \cite[Exercise 1.5.19.]{BH} we have $^*\Hom_S(M (-i), N(-j)) \cong {}^*\Hom_S(M, N)(i-j)$. Therefore $${}^*\Ext^i_S(M(-i),N(-j)) \cong {}^*\Ext^i_S(M,N)(i-j).$$

It can be easily shown that if $M$ is finitely presented, then ${}^*\Hom_S(M,N) = \Hom_S(M,N)$. If $S$ is both left and right Noetherian and $M$ is finitely generated, then 
\[{}^*\Ext^i_S(M,N)=\Ext^i_S(M,N),~~~~~~~~~~~~~~~\forall i \geq 0.\]

\s 
Let $M$ and $N$ be any two graded left $A_n(K)$-modules. We can define a left $A_n(K)$-module structure on ${}^*\Hom_R(M,N)$ extending the natural graded $R$-module structure by setting 
\[(\partial_j \cdot \varphi)(m) = \partial_j \cdot \varphi(m) - \varphi(\partial_j \cdot m).\]
Let $\varphi \in \Hom_i(M,N)$. Then for any $m \in M_n$ we have $\partial_j \cdot \varphi(m), \varphi(\partial_j \cdot m) \in M_{n+i-1}$ (as $\partial_j \cdot m \in M_{n-1}$). Thus $(\partial_j \cdot \varphi)(M_{n}) \subseteq M_{n+i-1}$ and hence $\partial_j \cdot \varphi \in \Hom_{i-1}(M,N).$ 

\begin{lemma}\label{inj-D-R}
Let $E$ be a ${}^*$injective $A_n(K)$-module, then $E$ is a ${}^*$injective $R$-module.
\end{lemma}

\begin{proof}
We know that
\begin{align*}
{}^*\Hom_R(-,E)&={}^*\Hom_R(-,{}^*\Hom_{A_n(K)}(A_n(K),E)), \\
&\cong {}^*\Hom_{A_n(K)}(A_n(K) \otimes_R-,E).
\end{align*}
Since $A_n(K)$ is free as graded right $R$-module, $A_n(K)\otimes_R -$ is exact. So ${}^*\Hom_{A_n(K)}(A_n(K)\otimes_R-,E)$ is exact and hence ${}^*\Hom_R(-,E)$ is exact. It follows that $E$ is a ${}^*$injective $R$-module.
\end{proof}

Consider a graded injective resolution of $N$ as an $A_n(K)$-module
\[\mathbb{I}^{\bullet}:~ 0 \to N \to I_0 \to I_1 \to \cdots \to I_r \to \cdots .\]
In view of Lemma \ref{inj-D-R}, $\mathbb{I}^{\bullet}$ is also a graded injective resolution of $N$ as an $R$-module. Therefore ${}^*\Ext_R^i(M,N) \cong H^i({}^*\Hom_R(M,\mathbb{I}^{\bullet}))$ have a natural structure of graded $A_n(K)$-modules.

\s Note that $M\otimes_R N$ is an $A_n(K)$-module with following action of $\partial_i$:  \[
\partial_i(m\otimes n) = (\partial_i m)\otimes n + m \otimes (\partial_i n).\]
Now consider  a free resolution  of the left $A_n(K)$-module $M$
\[\mathbb{F}_{\bullet} \colon~   \cdots \rt F_n \rt \cdots \rt F_1 \rt F_o \rt M \rt 0,\]
where $F_i$ are free $A_n(K)$-modules. Since $A_n(K)$ is a free $R$-module,  the complex $\mathbb{F}_{\bullet}$ is a resolution of $M$ by free $R$-modules.
It follows that $\Tor^R_\nu(M, N) = H_\nu(\mathbb{F}_{\bullet}\otimes_R N)$ have a natural structure of  $A_n(K)$-modules (see \cite[p.\ 18]{BJ} ).


\section{Strongly generalized Eulerian Modules}

In this section we 
discuss properties of the new class of graded $A_n(K)$-modules called strongly generalized Eulerian (Definition \ref{SGE}).
Let $\mathcal{T}$ be a graded Lyubeznik functor on ${}^*\Mod(R)$. Applying our theory of strongly generalized Eulerian module we prove that $\mathcal{T}(M)$ is strongly generalized Eulerian for any strongly generalized Eulerian module $M$.

The following  properties of strongly generalized Eulerian modules have already been proved for the case of generalized Eulerian modules in \cite{TP6}. Since proofs are almost similar to that of Proposition 2.1, Proposition 3.2, Proposition 5.3 and Theorem 6.3 in \cite{TP6}, so here we  skip  proofs.

\begin{proposition}\label{extension}
Let $ 0 \to N \to M \to L \to 0$ be a short exact sequence of graded $A_n(K)$-modules. Then $M$ is strongly generalized Eulerian if and only if $N, L$ are strongly generalized Eulerian.
\end{proposition}

\begin{proposition}\label{concen}
Let $M$ be a strongly generalized Eulerian $A_n(K)$-module. Then the $A_{n-1}(K)$-module $H_i(\partial_n, M)(-1)$ is strongly generalized Eulerian for $i=0,1$.
\end{proposition}

\begin{proposition}\label{kosz}
Let $M$ be a strongly generalized Eulerian $A_n(K)$-module. Then  $H_i(X_n, M)$ is strongly generalized Eulerian $A_{n-1}(K)$-module for $i=0,1$.
\end{proposition}

\begin{proposition}\label{Tor(MN)} Let $M$ and $N$ are graded strongly generalized Eulerian $A_n(K)$-modules. Then $\Tor^R_{i}(M,N)$ is a strongly generalized Eulerian $A_n(K)$-module for all $i \geq 0$.
\end{proposition}

Now we show that for any finitely generated $A_n(K)$- module , the notions of generalized Eulerian and strongly generalized Eulerian are equivalent. For this, first we prove the following lemma:

\begin{lemma}\label{obs} Let $A_n(K)$ be the $n^{th}$ Weyl algebra. Then for any $e \in \Z$ and $i=1, \ldots, n$, we have 
\[(\e_n-e) \cdot X_i=X_i (\e_n-(e-1)) \quad \mbox{and} \quad (\e_n-e) \cdot \partial_i=\partial_i (\e_n-(e+1)).\]
\end{lemma}

\begin{proof}
Recall that $\partial_iX_i=1+X_i\partial_i$ and $\partial_jX_i=X_i \partial_j$ for all $j \neq i$. Thus we have $(X_i\partial_i)\cdot X_i=X_i(1+X_i\partial_i)$ and hence
\begin{align*}	
(\e_n-e) \cdot X_i&=\left(X_i\partial_i+ \sum_{\substack{j=1\\j \neq i}}^nX_j\partial_j-e \right) \cdot X_i,\\
&=X_i\left(X_i\partial_i+\sum_{\substack{j=1\\j \neq i}}^n X_j\partial_j-e+1 \right),\\
&=X_i(\e_n-(e-1)).
\end{align*}

Again $(X_i\partial_i)\cdot \partial_i=(\partial_iX_i-1)\cdot \partial_i=\partial_iX_i\partial_i-\partial_i=\partial_i(X_i\partial_i-1).$ Therefore
\begin{align*}	
(\e_n-e) \cdot \partial_i&=\left(X_i\partial_i+\sum_{\substack{j=1\\j \neq i}}^nX_j\partial_j-e \right) \cdot \partial_i,\\
&=\partial_i\left(X_i\partial_i+\sum_{\substack{j=1\\j \neq i}}^nX_i\partial_i-e-1\right),\\
&=\partial_i(\e_n-(e+1)).\end{align*}
\end{proof}

\begin{remark}\label{rmk-obs}
From Lemma \ref{obs}, we get the following:
\begin{align*}
&(i)~(\e_n-e) \cdot X_i\partial_j=X_i(\e_n-(e-1)) \cdot \partial_j=X_i\partial_j(\e_n-(e-1+1))=X_i\partial_j(\e_n-e),\\
&(ii)~(\e_n-e) \cdot \partial_jX_i=\partial_j(\e_n-(e+1)) \cdot X_i=\partial_j X_i(\e_n-(e+1-1))=\partial_jX_i(\e_n-e),\\
&(iii)~(\e_n-e) \cdot X_i^t=X_i(\e_n-(e-1)) \cdot X_i^{t-1}= \cdots= X_i^t(\e_n-(e-t)),\\
&(iv)~ (\e_n-e) \cdot \partial_j^t=\partial_j(\e_n-(e+1)) \cdot \partial_j^{t-1}= \cdots= \partial_j^t(\e_n-(e+t)),
\end{align*}
for all $1 \leq i, j \leq n$ and $t \geq 1$.
\end{remark}

\begin{proposition}\label{ge-equiv-sge}
Let $M$ be a finitely generated graded $A_n(K)$-module. Then $M$ is generalized Eulerian if and only if it is strongly generalized Eulerian. 
\end{proposition}

\begin{proof}
Converse holds from the definition. We only need to prove the direct implication. Let $\{m_1, \ldots, m_r\}$ be a finite set of homogeneous generators of $M$. Let $\deg m_i=e_i$ and $(\e_n-e_i)^{b_i} \cdot m_i=0$ for some $b_i>0$. Set $b=\max\{b_1, \ldots, b_r\}$. We want to show that for any $e \in \Z$, $(\e_n-e)^{b} \cdot m=0$ for all $m \in M_{e}$. 

Denote $\alpha=(\alpha_1, \ldots, \alpha_n) \in \N^n$ and $|\alpha|=\alpha_1+\cdots+\alpha_n$. Take $d=\sum_{|\alpha-\beta|=t}c_{\alpha,\beta} X^\alpha\partial^\beta \in A_n(K)$ where $X^\alpha\partial^\beta=X_1^{\alpha_1} \cdots X_n^{\alpha_n} \partial_1^{\beta_1} \cdots \partial_n^{\beta_n}$ and $c_{\alpha,\beta} \in K$. By Remark \ref{rmk-obs}, for any $a \geq 1$ we get that 

\begin{align*}
(\e_n-e)^a d &= \sum_{|\alpha-\beta|=t}c_{\alpha,\beta}(\e_n-e)^a \cdot X^\alpha\partial^\beta,\\
&= \sum_{|\alpha-\beta|=t}c_{\alpha,\beta} X^\alpha(\e_n-(e-|\alpha|))^a\cdot \partial^\beta,\\
&= \sum_{|\alpha-\beta|=t}c_{\alpha,\beta} X^\alpha\partial^\beta(\e_n-(e-|\alpha|+|\beta|))^a,\\
&= \sum_{|\alpha-\beta|=t}c_{\alpha,\beta} X^\alpha\partial^\beta(\e_n-(e-t))^a,\\
&=d (\e_n-(e-t))^a.
\end{align*}
If $m=\sum_{i=1}^r d_im_i \in M_e$, then $\deg d_i=\deg m-\deg m_i=e-e_i$. Now we have
	\begin{align*}
	(\e_n-e)^{b} \cdot m &= \sum_{i=1}^r (\e_n-e)^{b} \cdot d_im_i,\\
	&=\sum_{i=1}^r d_i(\e_n-(e-(e-e_i)))^{b} \cdot m_i,\\
	&=\sum_{i=1}^r d_i (\e_n-e_i)^{b-b_i}(\e_n-e_i)^{b_i} \cdot m_i,\\
	&=0.
	\end{align*}
The result follows.
\end{proof}


\begin{remarks}\label{genarators-stg-gen}

\noindent
\begin{itemize}\item[(i)] Recall that any holonomic $A_n(K)$-module is finitely generated as an $A_n(K)$-module. So Lemma \ref{ge-equiv-sge} implies that for these modules both the notions generalized Eulerian and strongly generalized Eulerian are equivalent. 

\item[(ii)] If $M$ is a finitely generated graded $A_n(K)$-module with a finite generating set $\{m_1, \ldots, m_r\}$, then by Lemma \ref{ge-equiv-sge} we get that  $M$ is strongly generalized Eulerian if and only if $m_i$ satisfies \eqref{strong-ge-def} for every $i$. In particular, $A_n(K)/J$ is strongly generalized Eulerian for a homogeneous left ideal $J$ in $A_n(K)$ if and only if $\overline{1} \in A_n(K)/J$ satisfies \eqref{strong-ge-def} if and only if $\e_n^a \cdot 1=\overline{ \e_n^a}=0 \in A_n(K)/J$, i.e., $\e_n^a \in J$ for some $a>0$.  
\end{itemize}
\end{remarks}

From \cite[Theorem 1.7]{TP5} we have $\TT(R)$ is generalized Eulerian $A_n(K)$-module. 
Now we prove the following:

\begin{theorem}\label{lyu-sg}
Let $\TT$ be a graded Lyubeznik functor on ${}^*\Mod(R)$. Let $M$ be a strongly generalized Eulerian $A_n(K)$- module. Then $\TT(M)$ is strongly generalized Eulerian $A_n(K)$-module.
\end{theorem}

\begin{proof}
First we claim that if $M$ is strongly generalized Eulerian, then $M_f$ is strongly generalized Eulerian for any homogeneous polynomial $f \in R$. Let $m \in M$ be any homogeneous element. Since $M$ is strongly generalized Eulerian, $(\e_n-|m|)^am=0$ for some $a>0$. Notice ${\ensuremath \mid \frac{m}{f^i}\mid}=|m|-|f^i|$. By Lemma \ref{loca}, we get that 
\[(\e_n-|m|+|f^i|)^a\frac{m}{f^i}=\frac{1}{f^i}(\e_n-|m|)^am=0.\] 
Thus $M_f$ is strongly generalized Eulerian for any homogeneous polynomial $f \in R$. This follows that $M_{f_1\ldots f_s}$ is also strongly generalized Eulerian for homogeneous polynomials ${f_1,\ldots, f_s}$  in $R$. Let $I$ be an ideal of $R$ generated by ${f_1,\ldots, f_s}\in R$ 
and $C_\bullet$ be the $\check{C}$ech compex of $M$ with respect to $I$.  Then by the $\check{C}$ech complex characterization
of local cohomology and Proposition \ref{extension}, $H^i_{I}(M)$ is strongly generalized Eulerian $A_n(K)$-module for all $i\geq 0$. 

Let $Y=Y''-Y'$, where $Y', Y''$ are homogeneous closed sets of Spec$(R)$. Then $H^i_{Y'}(M)$ and $H^i_{Y''}(M)$ are both strongly generalized Eulerian $A_n(K)$-module. It follows from Proposition \ref{extension} that  $H^i_Y(M)$ is a strongly generalized Eulerian $A_n(K)$-module. Since $\TT=\TT_1\circ \TT_2\circ \ldots \circ \TT_s$ and each $\TT_i(M)$ is a graded $A_n(K)$-submodule of $H^i_Y(M)$, by Proposition \ref{extension} and using induction it follows that $\TT(M)$ is a strongly generalized Eulerian $A_n(K)$-module. 
\end{proof}

\begin{corollary}
Let $\TT$ be a graded Lyubeznik functor on ${}^*\Mod(R)$. Then $\TT(R)$ is strongly generalized Eulerian $A_n(K)$-module.

\end{corollary}

\begin{proof}
This follows immediately  from Theorem \ref{lyu-sg} as $R$ is  Eulerian $A_n(K)$-module .
\end{proof}


\section{Matlis dual of Graded Module \texorpdfstring{$M$}{M}}
We know that all the de Rham cohomology, Koszul homology and local cohomology modules are generalized Eulerian and behave nicely. In this section, we wish to investigate similar properties for $M^{\vee}$, the graded Matlis dual of $M$ using the notion of strongly generalized Eulerian modules.

\begin{definition} Let $R=K[X_1, \ldots, X_n]$ be a polynomial ring in $n$ variables. Let $\m=(X_1, \ldots, X_n)$ be the graded unique homogeneous maximal ideal of $R$. Then the Matlis dual of a graded $R$-module $M$ is defined by \[M^{\vee} ={}^* \Hom _K(M, K).\] 
\end{definition}

We can view $M$ as graded $K$-module and then $M^\vee$ is a graded $K$-module with grading \[(M^{\vee})_i=\Hom_K(M_{-i}, K)\] for all $i \in \Z$, see \cite[3.6, p. 141]{BH}. Clearly $M^\vee$ has a natural structure as a graded $R$-module. We can consider $K \cong R/\m$ as graded $R$-module with $K_0=K$ and $K_i=0$ for all $i \neq 0$. Let $^*E_R(K)$ denote the graded injective hull of $K$. By \cite[Proposition 3.6.16]{BH} we have that
\begin{align*}
& R^\vee= {}^* \Hom _K(R, K)\cong {}^*E_R(K),\\
\text{and} \quad & M^\vee={}^* \Hom _K(M, K)\cong{}^* \Hom _R(M, R^\vee)\cong{}^* \Hom _R(M, {}^*E_R(K))
\end{align*}
for all graded $R$-modules $M$, where the map \[\varphi: {}^* \Hom _R(M, {}^* \Hom_K(R, K)) \to {}^* \Hom _K(M, K)\] is defined by $\varphi(\alpha)(x)= \alpha(x)(1_R)$ for all $\alpha \in {}^* \Hom _R(M, {}^* \Hom _K(R, K))$ and $x \in M$.

\s \label{trnsp}{\bf Transposition of $A_n(K)$:}

The standard transposition of $A_n(K)$ is $\tau: A_n(K) \to A_n(K)$ defined by \[\tau(h\partial^{\alpha})=(-1)^{|\alpha|} \partial^{\alpha} h\] for all $h \in R=K[X_1, \ldots, X_n]$ and $\alpha \in \N^n$. 
Note that 
	\begin{itemize}
	\item[(i)]
	$\tau(ab)=\tau(b)\tau(a)$ for all $a,b \in A_n(K)$,
	\item[ (ii)] $\tau^2=\id_{A_n(K)}$.
\end{itemize} 
From this we get that	
\begin{enumerate}[\rm (1)]
\item $\tau(a^n)=\tau(a)^n$ for any $n \in \N$ and $a \in A_n(K),$
	
\item $\tau(\e_n)=\tau(\sum_{i=1}^n X_i \partial_i)= \sum_{i=1}^n \tau(X_i \partial_i)= \sum_{i=1}^n -\partial_iX_i= -\sum_{i=1}^n (X_i\partial_i+1)= -\e_n-n$ (as $\partial_i X_i- X_i\partial_i=1$ for all $i$).
\end{enumerate}

\s{\bf Transposition of left $A_n(K)$-module:}\label{trans}

For any left $A_n(K)$-module $M$, the {\it transpose} $M^{\sharp}$ of $M$ is a right $A_n(K)$-module such that 
\begin{enumerate}[\rm (i)]
\item $M^\sharp=M$ as abelian groups,
\item the right $A_n(K)$-action $\star$ on $M^{\sharp}$ is defined as $m \star a= \tau(a) \cdot m$ for all $a \in A_n(K)$ and $m \in M^{\sharp}$. 
\end{enumerate}
If $M=\bigoplus_{i \in \Z} M_i$ is a graded left $A_n(K)$-module, then $M^{\sharp}$ is also graded with $M^{\sharp}=\bigoplus_{i \in \Z} M_i$.
\begin{remark}
For any graded $A_n(K)$-module $M$, one can give a left $A_n(K)$-module structure on the module ${}^* \Hom _K(M, K)$ extending the natural $R$-module structure by 
\[(\partial_i \cdot \varphi)(m)= \varphi (\tau(\partial_i) \cdot m)\] 
for all $m \in M$ and $\varphi \in {}^* \Hom _K(M, K)$, as in \cite[Section 3, p. 9]{SZ}. Note that if $\varphi_1=\varphi_2$, then $(\partial_i \cdot \varphi_1-\partial_i \cdot \varphi_2)(m)=(\partial_i \cdot (\varphi_1-\varphi_2))(m)=(\varphi_1-\varphi_2)(\tau(\partial_i) \cdot m)=\varphi_1 (\tau(\partial_i) \cdot m)-\varphi_2 (\tau(\partial_i) \cdot m)=0$ (as $\varphi_1=\varphi_2$) for all $m \in M$ and hence $\partial_i \cdot \varphi_1=\partial_i \cdot \varphi_2$. Therefore the action is well-defined. 
\end{remark}

Now we prove the following result which is inspired by Proposition 3.9 in \cite{SZ}:

\begin{lemma}\label{matlis}
If $M$ is a strongly generalized Eulerian $A_n(K)$-module, then so is $M(-n)^{\vee}$.
\end{lemma}

\begin{proof}
Since $M$ is a strongly generalized Eulerian, for $m\in M(-n)_{-i}=M_{-n-i}$ there exists $a >0$ such that $(\e_n-(-n-i))^{a}\cdot m=0$. Thus for each $\varphi \in (M(-n)^{\vee})_i=\Hom_K(M(-n)_{-i}, K)$ and each $m \in M(-n)_{-i}$, we have 
\begin{align*}
	((\e_n-i)^{a}\cdot \varphi)(m)&= \varphi(\tau((\e_n-i)^{a}) \cdot m)\\
	&=\varphi((\tau(\e_n-i))^{a} \cdot m) \quad \text{as } \tau(a^n)=\tau(a)^n\\
	&=\varphi((\tau(\e_n)-i)^{a} \cdot m)\\
	&=\varphi((-\e_n-n-i)^{a} \cdot m)\\
	&=(-1)^{a}\varphi((\e_n+n+i)^{a} \cdot m)=0.
	\end{align*}
	Hence $(\e_n-i)^{a}\cdot \varphi=0$. The result follows. 
\end{proof}

The following result is well-known.
\begin{lemma}
Let $M$ and $N$ be graded $R$-modules. Then \[{}^*\Ext^i_R(M,N^\vee)= \Tor^R_i(M, N)^\vee \quad \mbox{for all $i \geq 0$}.\]
\end{lemma}

\begin{theorem}\label{Ext-R}
Let $N$ be strongly generalized Eulerian and $M$ be generalized Eulerian $A_n(K)$-modules. Then ${}^*\Ext^i_R(M,N)$ is generalized Eulerian for all $i \geq 0$..
\end{theorem}
\begin{proof} Since $N$ is a strongly generalized Eulerian $A_n(K)$-modules, so is $N(-n)^{\vee}$ by Lemma \ref{matlis}. In view of \cite[Theorem 6.3]{TP5} we get that $\Tor_i^R(M,N(-n)^{\vee})$ is generalized Eulerian. From Subsection \ref{ext-pre}
we have ${}^*\Ext^i_R(M(-i),N(-j))={}^*\Ext^i_R(M,N)(i-j)$. Hence 
\begin{align*}
\left[\Tor_i^R\left(M,N(-n)^{\vee}\right)\right] (-n)^{\vee}&=\left[\Tor_i^R(M(-n),N(-n)^{\vee})\right]^{\vee},\\
&={}^*\Ext^i_R(M(-n),N(-n)^\vee{}^\vee),\\
&={}^*\Ext^i_R(M(-n),N(-n)),\\ 
&={}^*\Ext^i_R(M,N)
\end{align*}
is generalized Eulerian.
\end{proof}

\begin{theorem}\label{Ext(MN)} 
Let $M, N$ be strongly generalized Eulerian $A_n(K)$-modules. Then $^*\Ext^R_{i}(M,N)$ is a strongly generalized Eulerian $A_n(K)$-module for all $i \geq 0$.
\end{theorem}
\begin{proof}
From Proposition \ref{Tor(MN)} and Lemma \ref{matlis} we get that $\Tor_i^R(M,N(-n)^{\vee})$ is strongly generalized Eulerian. Following the same method as Theorem \ref{Ext-R} one can easily show that $^*\Ext^R_{i}(M,N)$ is strongly generalized Eulerian.
\end{proof}

\newpage
\begin{remarks}
\noindent    
\begin{enumerate}[\rm (i)]
\item If $M,N$ are holonomic $A_n(K)$-modules, then $\Tor_i^R(M,N)$ is holonomic for all $i \geq 0$ by \cite[Theorem 1.6.4]{BJ}. In view of \cite[Example 3.14]{SZ} we get that $M(-n)^{\vee}$ need not be holonomic even if $M$ is holonomic. So from the proof of Theorem \ref{Ext-R} we can't say whether $\Ext^i_R(M,N)$ is holonomic or not.

\item If $\Ext^i_R(M,N)$ is supported at $\m$, then $\Ext^i_R(M,N) \cong E_R(K)(n)^\nu$ for some $\nu$ (possibly infinite), see\cite[Proposition 2.3, Theorem 2.4]{Lyu-1}. In particular, if $M=H^i_I(R)$ and $N=H_J^j(R)$ for some homogeneous ideals $I,J$ in $R$ and $\sqrt{I+J}=\m$, then there exists some $\nu_i \geq 0$ (possibly infinite) such that $\Ext^i_R(M,N) \cong E_R(K)(n)^{\nu_i}$.
\end{enumerate}
\end{remarks}


\section{\texorpdfstring{$\Ext$}{Ext} functor of \texorpdfstring{$A_n(K)$}{An(K)}-modules}
In this section, we provide a proof of our main result Theorem \ref{conj} which has already been conjectured  by  the first  and  third author in \cite{TP5} with strong evidence. To prove the result,  first we show the existence of  a first quadrant graded spectral sequence which is nontrivial and plays a vital role in the proof. 

The following result from \cite{PS} plays a very crucial role to prove  Theorem \ref{conj}:
\begin{lemma}\cite[p. 11]{PS}\label{ps-lem}
Let $M_2, M_3$ be left $A_n(K)$-modules and $M_1$ be a right $A_n(K)$-module. Then
\[\Hom_{A_n(K)}(M_1 \otimes_{R} M_2, M_3) \cong \Hom_{A_n(K)}(M_2, \Hom_R(M_1, M_3)).\]
\end{lemma}

\begin{remark}
Similarly one can easily show the  graded version of the above result:
If $M_2, M_3$ are graded left $A_n(K)$-modules and $M_1$ is a graded right $A_n(K)$-module, then
\[{}^*\Hom_{A_n(K)}(M_1 \otimes_{R} M_2, M_3) \cong {}^*\Hom_{A_n(K)}(M_2, {}^*\Hom_R(M_1, M_3)).\]
\end{remark}

Now we prove the following theorem which is the key result of this section and helps us to give the proof of the Theorem \ref{conj}: 

\begin{theorem}\label{conv}
Let $M,N$ and $L$ be graded $A_n(K)$-modules. Then there exists a first quadrant graded spectral sequence
\[{}^*\Ext^{p}_{A_n(K)}(R,{}^*\Ext^{q}_{R}(M,N)) \implies {}^*\Ext^{p+q}_{A_n(K)}(M,N).\]
\end{theorem}

\begin{proof}
Let $\mathbb{P}_{\bullet}$ be a finitely generated graded free resolution of $R$ as a $A_n(K)$-module and $\mathbb{I}^{\bullet}$ be a graded resolution of $N$ by ${}^*$injective $A_n(K)$-modules. Note that $\mathbb{I}^{\bullet}$ is also a graded injective resolution of $N$ as an $R$-module by Lemma \ref{inj-D-R}. From Lemma \ref{ps-lem} we have a graded isomorphism of bi-complexes 
\[{}^*\Hom_{A_n(K)}(\mathbb{P}_{\bullet},{}^*\Hom_R(M,\mathbb{I}^{\bullet})) \cong {}^*\Hom_{A_n(K)}(\mathbb{P}_{\bullet} \otimes_R M, \mathbb{I}^{\bullet}).\]
Set $X:={}^*\Hom_{A_n(K)}(\mathbb{P}_{\bullet},{}^*\Hom_R(M,\mathbb{I}^{\bullet})$ and $Y:={}^*\Hom_{A_n(K)}(\mathbb{P}_{\bullet} \otimes_R M, \mathbb{I}^{\bullet}).$ Since $X \cong Y$ as bigraded complexes, so we have $\Tot(X) \cong \Tot(Y)$. 

Now we compute the cohomology of $\Tot(Y)$. Set $E^{pq}={}^*\Hom_{A_n(K)}(\mathbb{P}_p \otimes_R M, \mathbb{I}^q)$.  Also we have 
\[
H_p(\mathbb{P}_p \otimes_R M) = \Tor_p^R(R,M)=
\left\{
\begin{array}{ll}
0  & \mbox{if } p > 0 \\
M  & \mbox{if } p=0
\end{array}
\right.\]
Now $\mathbb{I}^q$ is a ${}^*$injective $\D$-module. So
\[H^p({}^*\Hom_{A_n(K)}(\mathbb{P}_{\bullet} \otimes_R M, \mathbb{I}^q))=
\left\{
\begin{array}{ll}
{}^*\Hom_{A_n(K)}(M,\mathbb{I}^q) & \mbox{if } p = 0 \\
0  & \mbox{if } p \neq 0
\end{array}
\right.\]
It follows that 
\[{}^{I}E^{pq}_1=
\left\{\begin{array}{ll}
{}^*\Hom_{A_n(K)}(M,\mathbb{I}^q) & \mbox{if } p = 0 \\
0  & \mbox{if } p \neq 0
\end{array}
\right.\]
Now we take cohomology (arrows columns) and get 
\[{}^{I}E^{pq}_2=
\left\{\begin{array}{ll}
{}^*\Ext^q_{A_n(K)}(M,N) & \mbox{if } p = 0 \\
0  & \mbox{if } p \neq 0
\end{array}
\right.\]
The spectral sequence degenerates and hence 
$H^n(Y)={}^*\Ext^n_{A_n(K)}(M,N).$

Now we compute the cohomology of $\Tot(X)$. Set $E^{pq}={}^*\Hom_{A_n(K)}(\mathbb{P}_{q},{}^*\Hom_R(M,\mathbb{I}_{p}))$. Taking cohomology (arrows rows), we get
$H^p({}^*\Hom_R(M, \mathbb{I}^{\bullet}))={}^*\Ext^p_R(M,N).$ Since $\mathbb{P}_q$ is a finitely generated free $A_n(K)$-module, so we have 
$$H^p({}^*\Hom_{A_n(K)}(\mathbb{P}_{q},{}^*\Hom_R(M,\mathbb{I}^{\bullet}))) \cong {}^*\Hom_{A_n(K)}(\mathbb{P}_{q},{}^*\Ext^p_R(M,N)).$$
Hence 
\[{}^{I}E^{pq}_1={}^*\Hom_{A_n(K)}(\mathbb{P}_{q},{}^*\Ext^p_R(M,N)).\]
By taking cohomology (arrows columns), we get 
${}^{I}E^{pq}_2={}^*\Ext^q_{A_n(K)}(R,{}^*\Ext^p_R(M,N)).$
It follows that 
\[{}^*\Ext^{p}_{A_n(K)}(R,{}^*\Ext^{q}_{R}(M,N)) \implies {}^*\Ext^{p+q}_{A_n(K)}(M,N).\]
\end{proof}

\begin{lemma}\label{rel-ext-kos}
Let $L$ be a generalized Eulerian $A_n(K)$-module {\rm(}not necessarily holonomic{\rm)}. Then \[\Ext^{i}_{A_n(K)}({}^lR,L) \cong H^i(\partial, L)(-m).\]
\end{lemma}
\begin{proof}
	Note that $H^i(\partial, L) \cong \Tor^{A_n(K)}_{n-i}(R^r,L)$ is an isomorphism as graded $K$-vector spaces. Moreover, by \cite[Chapter 2, 7.5]{BJ} we have a graded isomorphism \[\Ext^i_{A_n(K)}({}^l R,L) \cong \Tor^{A_n(K)}_{n-i}(({}^lR)^{\dag},L),\] where $(-)^{\dag}=\Ext^n_{A_n(K)}(-,A_n(K))$. By \cite[Lemma 8.8]{TP5} we have $({}^lR)^{\dag} \cong R^r(-n)$. Therefore 
	\begin{align*}
	\Ext^i_{A_n(K)}({}^l R,L) &\cong \Tor^{A_n(K)}_{n-i}(R^r(-n),L), \\
	&\cong \Tor^{A_n(K)}_{n-i}(R^r,L)(-n), \\
	&\cong H^i(\partial, L)(-n).\end{align*}
The result follows.\end{proof}

We are now in a position to give a proof of our main result: 

\begin{theorem}\label{conj1}
Let $M$ and $N$ be non-zero left holonomic generalized Eulerian $A_n(K)$-modules. Fix  $i \geq 0$. Then the graded $K$-vector space $\Ext^i_{A_n(K)}(N,M)$ is concentrated in degree zero. 
\end{theorem}

\begin{proof}
In view of Lemma \ref{rel-ext-kos}, we have \[\Ext^{p}_{A_n(K)}({}^lR,{}^*\Ext^{q}_{R}(M,N)) \cong H^p(\partial, {}^*\Ext^{q}_{R}(M,N))(-n).\] Since $N$ is holonomic, $N$ is strongly generalized Eulerian by Lemma \ref{ge-equiv-sge}. So by Theorem \ref{Ext-R}, we get that $\Ext^{q}_{R}(M,N))$ is generalized Eulerian. Therefore $H^p(\partial, \Ext^{q}_{R}(M,N))$ is concentrated in degree $-n$ by \cite[Theorem 3.1]{TP6}. 

Since from Theorem \ref{conv} we have a graded convergence 
\[{}^*\Ext^{p}_{A_n(K)}(R,{}^*\Ext^{q}_{R}(M,N)) \implies {}^*\Ext^{p+q}_{A_n(K)}(M,N),\] it follows that $\Ext^i_{A_n(K)}(M,N)$ is concentrated in degree zero.
\end{proof}

As a consequence, we have the following surprising result:

\begin{corollary}
If $M$ and $N$ are graded holonomic generalized Eulerian left $A_n(K)$-modules, then any extensions \[\eta_i: ~ 0 \to N(i) \to Y \to M \to 0\] splits for all $i \neq 0$.
\end{corollary}

\begin{proof}
This follows immediately  from Theorem \ref{conj1} as $\Ext^1_{A_n(K)}(M,N)_i=0$ for all $i \neq 0$. 
\end{proof}

\begin{remark}
From \cite[Proposition 2.2]{TP6} we know that if $N$ is a non-zero generalized Eulerian $A_n(K)$-module. Then the shifted module $N(i)$ is not a generalized Eulerian for $i \neq 0$.
\end{remark}


\section{Acknowledgments}
We would like to thank \href{https://webusers.imj-prg.fr/~pierre.schapira/}{Prof. Pierre Schapira}, Sorbonne University, Institute of Mathematics of Jussieu for some useful conversations.

The first author thanks SERB, Department of Science and Technology, Govt. of India for the project  grant MATRICS (Project No.
MTR/2017/000585).
The second author thanks UGC, Govt. of India for providing financial support for this study. The third author is partially supported by Early Career Research Award funded
by SERB, Department of Science and Technology, Govt. of India (Project No.
ECR/2017/000963).

\let\cleardoublepage
\end{document}